

\magnification=\magstep1
\hsize=6.3truein
\input amssym.def
\input amssym.tex
\def\scong{{\scriptstyle\|}\lower.2ex\hbox{$\wr$}}

\def\Br{\mathop{\rm Br}\nolimits}

\def\Cor{\mathop{\rm Cor}\nolimits}

\def\rtimes{\mathop{\times\!\!{\raise.2ex\hbox{$\scriptscriptstyle|$}}}
	\nolimits} 
\def\proof{\noindent{\it Proof.}\quad}
\def\blackbox{\hbox{\vrule width6pt height7pt depth1pt}} 
\outer\def\Demo #1. #2\par{\medbreak\noindent {\it#1.\enspace}
	{\rm#2}\par\ifdim\lastskip<\medskipamount\removelastskip
	\penalty55\medskip\fi}
\def\qed{~\hfill\blackbox\medskip}
\overfullrule=0pt
\def\Br{\mathop{\rm Br}\nolimits}
\def\Cor{\mathop{\rm Cor}\nolimits}

\def\hangbox to #1 #2{\vskip1pt\hangindent #1\noindent \hbox to #1{#2}$\!\!$}

\pageno=0
\footline{\ifnum\pageno=0\hfill\else\hss\tenrm\folio\hss\fi}
\topinsert\vskip1.8truecm\endinsert
\centerline{GENERIC ALGEBRAS WITH INVOLUTION OF DEGREE 8m}
\vskip6pt
$${\vbox{\halign{\hfil\hbox{#}\hfil\qquad&\hfil\hbox{#}\hfil\cr
$$David J. Saltman$^*$\cr
Department of Mathematics\cr
The University of Texas\cr
Austin, Texas 78712\cr}}}$$
\bigskip
\centerline{and}
\bigskip
$${\vbox{\halign{\hfil\hbox{#}\hfil\qquad&\hfil\hbox{#}\hfil\cr
$$Jean - Pierre Tignol$^*$\cr
Institut de Math\'ematique Pure et Appliqu\'ee\cr
Universit\'e Catholique de Louvain\cr
Louvain-la-Neuve, Belgium\cr}}}$$
\vskip16pt
\vskip3pt
{\narrower\smallskip\noindent
{\bf Abstract} The centers of the generic central simple algebras with 
involution are interesting objects in the theory of 
central simple algebras. These fields also arise 
as invariant fields for linear actions of projective 
orthogonal or symplectic groups. In this 
paper, we prove that when the characteristic is 
not 2, these fields are retract rational, 
in the case the degree is $8m$ and $m$ is odd. 
We achieve this by proving the equivalent lifting 
property for the class of central simple algebras 
of degree $8m$ with involution. A companion paper ([S3]) 
deals with the case of $m$, $2m$ and $4m$ 
where stronger rationality results are proven.
\bigskip

\medskip

\smallskip}

\footnote{}{*The first author is grateful for support under NSF grant 
DMS-9970213. The first author would also like to thank the 
Universit\'e Catholique de Louvain for its hospitality. The second 
author is supported in part by the National Fund for Scientific 
Research (Belgium) and by the TMR network 
``Algebraic K - Theory, Linear Algebraic Groups and Related Structures''
(ERB FMRX CT 97-0107)}
\vfill\eject 
In this paper $F$ will always be an infinite field of characteristic not 2. 
Let ${\cal G}$ be an algebraic group over $F$ and $V$ 
an algebraic $F$ representation, by which we mean 
there is an algebraic group morphism ${\cal G} \to GL_F(V)$. 
There is considerable interest in the structure, and 
more specifically in the 
rationality, of the invariant field $F(V)^{\cal G}$, 
where ${\cal G}$ has its natural action on the field of rational 
functions $F(V)$ of $V$. For specific groups and $V$, 
this question has particular significance. For example, 
consider ${\cal G} = PGL_n(F) = GL_n/F^*$ and 
$V = M_n(F) \oplus \ldots \oplus M_n(F)$ ($r$ times) where the action 
of $PGL_n(F)$ on $V$ is induced by diagonal conjugation. 
Then the invariant field $F(V)^{PGL_n}$ is the center of a 
generic division algebra $UD(F,n,r)$ (e.g. [LN] sec. 14). 

In $PGL_n$ there are subgroups and for some of these subgroups 
the corresponding invariant field is also of importance. 
We will be particularly interested in the projective orthogonal 
groups $PO_n$ and projective symplectic groups 
$PSp_n$ (for $n$ even). Since we do not assume $F$ 
is algebraically closed, let us be precise here. Let 
$O_n(F) \subset GL_n(F)$ be the group of orthogonal 
matrices. That is, $O_n(F)$ is the group of matrices 
where $AA^T = I$, where $T$ is the transpose. Let 
$Sp_n(F)$ be the group of symplectic matrices, 
that is the group of matrices where $AA^S = I$  
and $S$ is the standard symplectic involution. 
For our purposes we can then define $PO_n(F)$ and 
$PSp_n(F)$  to be the image of $O_n(F)$ and $Sp_n(F)$ 
in $PGL_n(F)$. Note that, with this choice, 
$PO_n(F)$ and $PSp_n(F)$ may not be the group of $F$ rational 
points of the corresponding algebraic group, because 
the quotient groups may have $F$ points not in the image 
of the group of $F$ points of $O_n$ or $Sp_n$. 
To remedy this one could replace $O_n$ and $Sp_n$ by 
$GO_n$ and $GSp_n$, the corresponding groups of similitudes 
(e.g. [K-T] p. 153). However, for our purposes 
none of this matters. Our definition of $PO_n(F)$ and 
$PSp_n(F)$ yield a Zariski dense set of points in the 
corresponding groups over the algebraic closure of $F$, 
and so the invariant rings and fields are the same 
no matter what definition we take. 

In, for example, [R1] p. 183 there is a definition of generic 
algebras $UD_t(F,n,r)$ and $UD_s(F,n,r)$ with involution 
of orthogonal respectively symplectic type.  
By [P] p. 377-378, 
$F(V)^{PO_n}$ is the center, $Z_t(F,n,r)$, of $UD_t(F,n,r)$ 
while $F(V)^{PSp_n}$ is the center, $Z_s(T,n,r)$, of $UD_s(F,n,r)$. 
Thus the invariant fields of $PO_n$ and $PSp_n$ play the role 
in the theory of central simple algebras with involution 
that the invariant field of $PGL_n$ plays in the theory 
of central simple algebras. In particular, these invariant fields 
are natural objects to consider.

Though the original question we asked was about 
rationality, there is a weaker property which is closely tied 
to properties of central simple algebras. 
We say a field extension $K/F$ is {\bf retract rational} 
if and only if the following holds. $K$ is the field of fractions $q(S)$ 
of an $F$ algebra domain $S$, and there is a localized polynomial 
ring $F[\vec x](1/s) = F[x_1, \ldots, x_n](1/s)$ with $F$ algebra maps 
$f: S \to F[\vec x](1/s)$ and $g: F[\vec x](1/s) \to S$ 
such that $g \circ f: S \to S$ is the identity. 

The basic properties of retract rational field extensions 
are developed in [S]. Let us note one here. Define $K,K'$ 
to be stably isomorphic (over $F$) if and only if 
the following holds. For some $a$,$b$, 
the fields $K(x_1,\ldots,x_a)$ and $K'(y_1,\ldots y_b)$ 
are isomorphic over $F$, where the $x$'s and $y$'s 
ate transcendence bases. 
It is shown in [S] that if $K,K'$ 
are stably isomorphic, and $K/F$ is retract rational, 
then $K'/F$ is retract rational. In particular, 
stably rational (i.e. stably isomorphic to a rational extension) 
implies retract rational (but not conversely). Because 
of the above fact, we will talk about the retract rationality 
of the stable isomorphism class of a field extension $K/F$. 

Let us break to explain a little notation. 
The statement $A/K$ is a central simple algebra 
of degree $n$ means that $A$ is a simple algebra of dimension 
$n^2$ over its center $K$. If we say $D/K$ is a division 
algebra, we also mean $K$ is its center. 
If $A/K$ is central simple, we will write 
$K(A)$ to mean the function field of the Severi Brauer variety 
of $A$. That is, $K(A)$ is the Amitsur generic splitting field 
of $A$. Finally, 
suppose $A/K$ and $A'/K'$ are central simple algebras 
and $K(x_1,\ldots,x_a) \cong K'(y_1,\ldots,y_b)$ as in the 
definition of stable isomorphism. If some such isomorphism 
extends to an isomorphism $A \otimes_K K(x_1,\ldots,x_a) 
\cong A' \otimes_{K'} K'(y_1,\ldots,y_b)$, we say 
$A/K$ and $A'/K'$ are stably isomorphic.    

As mentioned above, $F(V)^{PO_n}$ and $F(V)^{PSp_n}$ 
are the centers of the so called generic algebras with 
orthogonal respectively symplectic involution. In particular, 
these fields are centers for generic objects for the 
class of central simple algebras with orthogonal 
respectively symplectic involutions. It follows that these 
are also generic objects for the class of central simple 
algebras of order dividing 2 in the Brauer group. 
This last fact is reflected in the result from [BS] 
we are about to quote in Theorem 1, describing $F(V)^{PO_n}$ and 
$F(V)^{PSp_n}$ as extensions of $F(V)^{PGL_n}$. 
Furthermore, in Theorem 2, we will confront more precisely 
what it means to be a generic object for a class 
of central simple algebras. 

To state it the result from [BS] we need, 
let $r$  be the number of direct summands in $V$ and 
$UD(F,n,r)/Z(F,n,r)$ the generic division algebra of degree 
$n$ in $r$ variables. Abbreviate $UD/Z = UD(F,n,r)/Z(F,n,r)$. 
Let $B_o$ be the central simple 
algebra of degree $n(n+1)/2$ in the Brauer class 
of $UD \otimes_Z UD$ and $B_s$ the central simple algebra 
in the same class of degree $n(n-1)/2$. 
Note that $B_o$ is written $s^2UD$ and 
$B_s$ is written $\lambda^2UD$ in [K-T] p. 33.

\proclaim Theorem 1. 
For any $n$, $F(V)^{PO_n} = Z_t(F,n,r)= Z(B_o)$. If $n$ 
is even (so $PSp_n$ is defined), $F(V)^{PSp_n} = Z_s(F,n,r) = Z(B_s)$. 

Let $D'$ be the division algebra in the class of $UD \otimes_Z UD$. 
Then by e.g. [LN] p. 93, $Z(B_o)$ and $Z(B_s)$ are, when defined, 
rational over $Z(D')$. In particular, $Z(B_o)$ is isomorphic 
to a field rational over $Z(B_s)$. Thus, to save ink, 
we will frequently only discuss $F(V)^{PO_n} = Z_t(F,n,r)$ 
since the other field is equivalent. 

The goal of this note is a result on retract rationality, 
which we prove by relating retract rationality to a property 
of algebras. To this end, let ${\cal A}_{2,n}$ be the class of Azumaya 
algebras $A/R$ of degree $n$ where $R \supset F$ and 
$A \otimes_R A \cong M_t(R)$ for the appropriate $t$. 
Note that this is a linear class in the sense of [LN] p. 76. 
We say ${\cal A}_{2,n}$ has the lifting property ([LN] p. 77)
if and only if the following holds. 
Assume $T$ is a local commutative $F$ algebra with residue 
field $K$ and $A/K$ is in ${\cal A}_{2,n}$.  
Then there is an Azumaya $B/T \in {\cal A}_{2,n}$ 
with $B \otimes_T K \cong A$. 

Lifting is important because of Theorem 2 to follow. 
But before we state the result, we recall a few notions from 
[LN] section 11. $UD_t = UD_t(F,n,r)$ can be identified 
with $UD \otimes_Z Z(B_o)$ and the center of both 
these algebras can be identified with $F(V)^{PO_n}$. 
Suppose $A/S$ is an Azumaya 
such that $q(S) = F(V)^{PO_n}$ and 
$A \otimes_S F(V)^{PO_n} = UD_t$. 
If $B/R \in {\cal A}_{2,n}$, we say $\phi: S \to R$ 
realizes $B$ if and only if $B \cong A \otimes_{\phi} R$. 
Note that $\otimes_{\phi}$ means that we treat $R$ as an 
$S$ module via $\phi$. 

We say $UD_t/Z(B_o)$ represents ${\cal A}_{2,n}$ (see [LN] p. 76)
if and only if the following holds. 
There is an $A/S$ Azumaya such that $S$ is finitely generated 
as an $F$ algebra, $q(S) = Z(B_o)$, $A \otimes_S Z(B_o) \cong 
UD_t$, and further the following holds. 
Assume $0 \not= s \in S$ and $B/K \in {\cal A}_{2,n}$ 
with $K$ a field. Then there is a $\phi: S(1/s) \to K$ 
realizing $B/K$.  
Note that if $A/S$ is as above, and 
$S' \subset F(V)^{PO_n}$ satisfies $q(S') = F(V)^{PO_n}$, 
then for some $0 \not= s' \in S'$ and some $A'/S'(1/s')$, 
$A'/S'(1/s')$ satisfies the same property. This is why we can 
view ``representing'' as a property of the algebra $UD_t/Z(B_o) 
= UD_t/F(V)^{PO_n}$.  
Also, it is clear that if $UD_t/F(V)^{PO_n}$ is stably isomorphic 
to a  $A/K$, and $A/K$ represents ${\cal A}_{2,n}$, 
then so does $UD_t/F(V)^{PO_n}$. Thus we can talk of the stable 
isomorphism class of $UD_t/F(V)^{PO_n}$ as 
representing ${\cal A}_{2,n}$. 

Another idea we recall is called ``local projectivity'' 
in [S], or (a slight variant) property v) in [LN] p. 76. 
We will use the version of this property from [LN], but the name 
local projectivity from [S]. 
Let $A/S$ be such that $q(S) = F(V)^{PO_n}$ and 
$A \otimes_S F(V)^{PO_n} = UD_t$. 
Suppose $B'/T \in {\cal A}_{2,n}$ and $T$ is a local ring 
with residue field $K$. Set $B = B' \otimes_T K$. 
Then $A/S$ is locally projective if and only if 
for any such $B'/T$ etc., and any $\phi \to K$ realizing $B/K$, 
there is a $\phi': S \to T$  realizing $B'/T$ such that 
the composition $S \to T \to K$ is $\phi$. 
Note that if $A/S$ is locally projective then so is 
$A(1/s)/S(1/s)$ for any $0 \not= s \in S$. 
Thus once again it is fair to talk about $UD_t/F(V)^{PO_n}$ 
being locally projective. Also it is clear that the 
property of being locally projective is preserved by stable 
isomorphisms. Thus, once again, we can talk about the 
stable isomorphism class of $UD_t/F(V)^{PO_n}$ as being 
locally projective. 

In [S] and [LN sec. 11] a general framework is described 
along with a result connecting lifting properties 
with retract rationality. This framework applies here  
and so we can show:  

\proclaim Theorem 2. The stable isomorphism classes of  
$F(V)^{PO_n}/F = Z_t(F,n,r)/F$ or 
(when $n$ even) 
$F(V)^{PSp_n}/F = Z_s(F,n,r)/F$ are retract rational if and only if 
${\cal A}_{2,n}$ has the lifting property. 

\proof 
By [LN] p. 77 it is enough show that $UD_t/F(V)^{PO_n}$ represents  
${\cal A}_{2,n}$ and is locally projective. By the 
above observations, we can replace $Z(B_o) = 
F(V)^{PO_n}$ by $K = Z(UD \otimes_Z UD)$, and 
$UD_t$ by $D = UD_t \otimes_{Z(B_o)} K$, because 
$K/Z(B_o)$ is rational (e.g. [LN] p. 93). 
 
In [S1] was defined a generic central simple algebra 
$D'/K'$ of degree $n$ and order dividing $t$. 
In that paper $D'/K'$ was shown to represent the class of Azumaya 
algebras with the same property. 
In the case of $t = 2$, it follows from [S2] p. 344 
that $D'/K'$ is rational over $D/K$, 
and so $D/K$ represents ${\cal A}_{2,n}$. 

In [LN] p. 105 it is shown that $UD/Z$ 
is locally projective for the class of Azumaya 
algebras of degree $n$. Let $A'/S'$, $q(S') = Z$, 
be an Azumaya algebra that realizes this property. 
Define $S \supset S'$ to be the affine ring of 
an affine open subset of the Severi-Brauer scheme of 
$A' \otimes_{S'} A'$ (e.g. [V]) and set $A = A' \otimes_{S'} S$. 
Then $q(S) = Z(UD \otimes_Z UD) = K$ by the naturality 
of the Severi-Brauer scheme. Furthermore, clearly 
$A \otimes_S K = D$. 
We claim that using $A/S$ one sees that $D/K$ is locally projective. 

Suppose $B'/T$ is in ${\cal A}_{2,n}$, $T$ is local with residue 
field $K$, and $B = B' \otimes_T K$. Assume $\phi: S \to K$ 
realizes $B/K$. Since $A'/S'$ is locally projective, 
there is a partial lifting $\phi'': S' \to T$ which realizes 
$B'$. That is, the restriction $\phi|_{S'}: S' \to K$ 
can be factored into $S' \to T \to K$ where the first map 
is $\phi''$.  The full map $\phi$ can be factored 
into $S \to S \otimes_{\phi''} T \to K$. Note that by the 
naturality of the Severi Brauer scheme, $S \otimes_{\phi''} T$ 
is the affine ring of the corresponding open subset, call it $U$, 
of the Severi Brauer scheme of $B' \otimes_T B'$. 
Thus $\phi$ defines a 
$K$ point on the Severi-Brauer variety of $B' \otimes_T B'$ 
which can be identified with a $K$ point of the Severi Brauer variety of 
$B \otimes_K B$.  
There is a transitive action by $(B \otimes_K B)^*$ on these $K$ points, 
and $(B' \otimes_T B')^*$ maps onto $(B \otimes_K B)^*$. 
By assumption, there is a $T$ point on the Severi-Brauer scheme 
of $B' \otimes_T B'$. It follows that the  
$K$ point given by $\phi$ is the image of a $T$ point of the Severi-Brauer 
scheme of $B' \otimes_T B'$. Since $T$ is local, 
the closure of this $T$ point includes the $\phi$ given $K$ 
point, and so this $T$ point is also in $U$. 
That is, there is a morphism $S \otimes_{\phi''} T \to T$ 
and the composition $\phi': S \to S \otimes_{\phi''} T \to T$ 
is the required lift for $\phi$. 
This proves local projectivity and hence Theorem 2. 

It is clear how we will use Theorem 2, but before we do that 
let us make one final reduction. 

\proclaim Lemma 3. Let $n = 2^rm$ where $m$ is odd. 
Then ${\cal A}_{2,n}$ has the lifting property 
if ${\cal A}_{2,2^r}$ has the lifting property. 

\proof 
If $A/K$ is in ${\cal A}_{2,n}$, then $A = A_2 \otimes A_m$ 
where $A_2$ has degree $2^r$ and $A_m$ has degree $m$ 
(e.g. [LN] p. 35). Since $A$ has order 2 in the Brauer group, 
and $A_m$ has order dividing $m$, it follows that $A_m$ 
must be split. That is, $A \cong M_m(A_2)$. 
It is now obvious that if ${\cal A}_{2,2^r}$ has the lifting 
property then so does ${\cal A}_{2,n}$.~\qed

We remark that the converse is also true, but 
to prove this would take us too far afield. 
To outline the argument, if $B/T$ is an Azumaya 
algebra over a local ring, then $B \cong M_s(D)$ 
where $D$ has no nontrivial idempotents. 
Moreover, there is only one such $D$, up to isomorphism, 
in the Brauer class of $B$. With this, one can copy the 
usual proof over a field, and show that $B \cong 
B_1 \otimes_T \ldots \otimes_T B_s$ where all the 
$B_i$ have prime power degree. With this 
background, the converse is clear.  

We can now state:  

\proclaim Theorem 4. Suppose $F$ is an infinite field of 
characteristic not 2 
and $n = 8m$ where $m$ is odd. Then the stable isomorphism classes of 
$F(V)^{PSp_n}$ and $F(V)^{PO_n}$ are retract rational over $F$. 
Equivalently, the stable isomorphism classes of 
the centers $Z_t(F,n,r)$ and $Z_s(F,n,r)$ 
of the generic algebras with orthogonal respectively 
symplectic involution are retract rational over $F$. 

Before we prove Theorem 4, we begin with another lemma.  
Let $R$ be a commutative ring. If $b_i \in R$ are finitely many elements, 
define $R(b_1^{1/2}, \ldots, b_s^{1/2})$ to be 
$R[x_1, \ldots, x_s]/< x_i^2 - b_i | i = 1, \ldots s >$. 
Note that we make the above definition even if some of the 
$b_i$ are squares. In particular, if $R$ is a field, 
$R(a_1^{1/2}, \ldots, a_s^{1/2})$ may not be a field 
but is a direct sum of fields. We recall: 

\proclaim Lemma 5. Let $T$ be a local $F$ algebra  
with residue field $K$. Suppose $a_i \in K^*$ and 
$a_i' \in T$ are preimages. Then $S = T(a_1'^{1/2}, \ldots, a_s'^{1/2})$ 
is a semilocal $F$ algebra which, modulo its Jacobson radical, is isomorphic to $L = K(a_1^{1/2}, \ldots, a_s^{1/2})$. 
In particular, $S^*$ maps onto $L^*$. 
$S/T$ is Galois with Galois group we can identify with the 
Galois group of $L/K$. Call this group $G$. There is an isomorphism 
$H^2(G,S^*) \cong \Br(S/T)$. 

\proof 
Since the $a_i'$ are invertible, it is easy to see $S/T$ 
is Galois and since Galois extensions are closed under 
specialization, one can identify this Galois group 
with that of $L/K$. The Jacobson radical of $S$ must be 
${\cal M}S$ where ${\cal M}$ is the maximal ideal of 
$T$. Since $L$ is a direct sum of fields, $S$ is semilocal. 
Of course, semilocal local rings have trivial Picard group, 
so $H^2(G,S^*) \cong \Br(S/T)$ by, e.g., [LN] p.45.~\qed 

If $A'$ is any $T$ algebra, and $T$ has residue field 
$K$, then we say $A'$ is a lift of $A = A' \otimes_T K$. 
When $A/K$ is central simple, we will only call $A'$ 
a lift if $A'/T$ is Azumaya. When $A/K$ is a commutative 
Galois extension with Galois group $G$, we will only say $A'$ 
is a lift if $A'/T$ is Galois with group $G$. 
Thus among the results of Lemma 5 is that 
$T(a_1'^{1/2},\ldots,a_s'^{1/2})$ is a lift of 
$K(a_1^{1/2},\ldots,a_s^{1/2})$. 

Let us also recall that if $R$ is any commutative ring containing 1/2, 
and $a,b \in R^*$, then one can form the Azumaya 
quaternion algebra $(a,b)_R = 
R \oplus R\alpha \oplus R\beta \oplus R\alpha\beta$ where 
$\alpha^2 = a$, $\beta^2 = b$, and $\alpha\beta = -\beta\alpha$. 
As implied, $(a,b)_R$ is Azumaya over $R$ of rank 4 (i.e. degree 2) 
([LN] p. 49). 
By e.g. [LN] p. 34, $(a,b)$ defines an element of order 2 in the Brauer 
group of $R$. 
Furthermore, $(a,b)_R \cong (b,a)_R \cong (a,N_S(\gamma)b)_R$ 
where $\gamma \in R(a^{1/2})^*$, $S = R(a^{1/2})$, and 
$N_S: R(a^{1/2}) \to R$ 
is the norm. If $R$ is semilocal, then $(a,b) \cong (a,c)$ 
implies $bc$ is a norm from $R(a^{1/2})$ by Lemma 5. 

Let $a \in R^*$ with $R$ as above, 
and $S = R(a^{1/2})$. Then the corestriction 
$\Cor_{S/R}: \Br(S) \to \Br(R)$ is defined (e.g. [LN] p. 55) 
and satisfies all the usual properties. In particular, 
if $a \in R^*$ and $b \in S^*$, then 
$\Cor_{S/R}((a,b)_S)$ is Brauer equivalent to $(a,N_S(b))_R$ 
(e.g. [LN] p. 57). Furthermore, if $A/R$ is Azumaya, 
$\Cor_{S/R}(A \otimes_R S)$ 
is Brauer equivalent to $A \otimes_R A$. 
Let $\sigma$ generate the Galois group of $S/R$. 
That is, $\sigma(a^{1/2}) = -a^{1/2}$. Suppose $B/S$ 
is Azumaya and let $\sigma(B)$ be the $\sigma$ twist. 
That is, $\sigma(B) = B \otimes_{\sigma} S$. 
We finally have $\Cor_{S/R}(B) \otimes_R S$ is Brauer equivalent 
to $B \otimes_S \sigma(B)$. 

We are finally ready to turn to the proof of Theorem 4. 
Of course, by Theorem 2 and Lemma 3 it suffices to prove 
${\cal A}_{2,8}$ has the lifting property. To this end, 
suppose $T$ is a local $F$ algebra with residue field 
$K$, and $D/K$ is a central simple algebra of degree 
8 and order 2 in the Brauer group. We must show 
that there is an Azumaya $D'/T$ such that $D' \otimes_T K \cong D$ 
and $D' \otimes_T D'$ is isomorphic to matrices 
over $T$. Note that since $T$ is local, this is equivalent 
to saying $D'$ has order dividing 2 in the Brauer group. 

By [R], $D$ has a maximal subfield  
of the form $K(a_1^{1/2},a_2^{1/2},a_3^{1/2})$. 
The centralizer of $L = K(a_1^{1/2})$ in $D$ is a division 
algebra of degree 4 with involution. Thus by e.g. [LLT] Proposition 5.2, 
this centralizer has the form $B = (a_2,x_2)_L \otimes_L (a_3,x_3)_L$. 

The corestriction of $[B]$ is Brauer equivalent to $D \otimes_K D$ 
and so must be trivial. But this corestriction 
is $(a_2,N_L(x_2))_K \otimes_K (a_3,N_L(x_3))_K$. 
In other words, $(a_2,N_L(x_2))_K \cong (a_3,N_L(x_3))$. 
By [T] p.267 or [A] Lemma 1.7, there is a $y \in K^*$ such that 
$(a_2,N_L(x_2)) \cong (y,N_L(x_2)) \cong (y,N_L(x_3)) \cong 
(a_3,N_L(x_3))$. Set $L_i = K(N_L(x_i)^{1/2})$ for $i = 2,3$ 
and $L_{23} = K(N_L(x_2x_3)^{1/2})$. Then 
there are $\mu_i \in L_i^*$ and $\mu_{23} \in L_{23}^*$ 
such that $a_2y = N_{L_2}(\mu_2)$, $y = N_{L_{23}}(\mu_{23})$, 
and $a_3y = N_{L_3}(\mu_3)$. The idea of this proof 
is that we can lift $a_1$, then the $x_i$, then $y$, and then 
$a_2$, $a_3$ so that all these relations still hold. 
The key idea is that we use the relations to define the lifts.  

Choose $a_1' \in T^*$ a preimage of $a_1$. Set $S = T(a_1'^{1/2})$, 
so $S$ is a lift of $L$. 
Choose $x_i' \in S^ *$ preimages of the $x_i$. 
Of course, $N_S(x_i')$ is a preimage of $N_L(x_i)$. 
Set $S_i = T(N_S(x_i')^{1/2})$ and $S_{23} = T(N_S(x_2'x_3')^{1/2})$. 
Of course, the $S_i$ and $S_{23}$ are lifts of the $L_i$ and $L_{23}$ 
respectively. Choose $\mu_i' \in S_i^*$ and $\mu_{23}' \in S_{23}^*$ 
preimages of the $\mu_i$ and $\mu_{23}'$ respectively. 

Set $y' = N_{S_{23}}(\mu_{23}')$. Clearly $y' \in T^*$ 
is a preimage of $y$. For $i = 2,3$, set 
$a_i' = N_{S_i}(\mu_i')y'^{-1} \in T^*$. 
Clearly, the $a_i'$ are preimages of the $a_i$. 
Set $B' = (a_2',x_2')_S \otimes_S (a_3',x_3')_S$. 
Of course, $B'$ is a lift of $B$. 
The corestriction $\Cor_{S/T}(B')$ is Brauer equivalent 
to $(a_2',N_S(x_2'))_T \otimes_T (a_3',N_S(x_3'))_T$. 
But $(a_2',N_S(x_2')) \cong (y',N_S(x_2')) \cong (y',N_S(x_3')) 
\cong (a_3',N_S(x_3'))$. It follows that 
$\Cor_{S/T}(B')$ is trivial. Tensoring up to 
$S$, we have $B' \otimes_S\sigma(B')$ is trivial 
where $\sigma$ generates the Galois group of $S/T$. 
Of course this means $B'$ and $\sigma(B')$ are Brauer 
equivalent. Since $S$ is semilocal, using [D] we have that 
$B' \cong \sigma(B')$. Alternatively, we can make the 
following argument. 
Both $B'$ and $\sigma(B')$ are split by $V = S(a_2'^{1/2},a_3'^{1/2})$. 
More precisely, both $B'$ and $\sigma(B')$ are crossed products 
(e.g. [OS] p. 88-90)
with respect to $V/S$. By [LN] p. 45, the corresponding cocycles 
are cohomologous, and so $B' \cong \sigma(B')$. 

The isomorphism $B' \cong \sigma(B')$ can be equivalently 
expressed as the existence of  
an $\alpha: B' \cong B'$ 
such that $\alpha$ is $\sigma$ semilinear. 
Since $\alpha^2$ is an $S$ automorphism, and $S$ is semilocal, 
$\alpha^2$ is an inner automorphism given by, say, $c \in B'^*$ 
(e.g. [LN] p. 16). 

Form the algebra $A' = B' \oplus B'u$ where $ub = \alpha(b)u$ 
for all $b \in B'$ and $u^2 = c$. Using e.g. [LN] p. 12 
it is easy to see that $A'/T$ is Azumaya over $T$ of degree 8, 
and the centralizer, 
in $A'$, of $S \subset B'$ is $B'$. Thus (e.g. [LN] p. 24) 
$A'/T$ defines a preimage of $B'$ in the Brauer group 
of $T$. In particular, $A' \otimes_T A'$ is Brauer equivalent to 
$\Cor_{S/T}(B')$ and so $A'$ has order 2 in the Brauer group. 

If $A = A' \otimes_T K$, then $A$ and $D$ have equal 
images in the Brauer group of $L$. That is, 
$M_2(A) \cong D \otimes_K (a_1,d)$ for some $d \in K^*$. 
Let $d' \in T$ be a preimage of $d$  and set $A'' = 
A' \otimes_T (a_1',d')$. Of course, the Brauer class 
of $A''$ is a preimage of the Brauer class of $D$. 
$A''$ contains the subalgebra 
$S \otimes_T S$. Since $S/T$ is Galois, $S \otimes_T S$ 
contains an idempotent $e$ such that $e(S \otimes_T S) \cong S$. 
Viewing $e \in A''$, it is easy to see that $D' = eA''e$ 
is Azumaya over $T$ of degree 8 and so $D'$ is a lift of 
$D$.~\qed

\bigskip
\leftline{References}
\def\hangbox to #1 #2{\vskip1pt\hangindent #1\noindent \hbox to #1{#2}$\!\!$}
\def\refn#1{\hangbox to 30pt {#1\hfill}}

\medskip
\refn{[A]}
Arason, J.K., 
{\it Cohomologische Invarianten Quadratischen Formen}, 
J. Algebra {\bf 36} (1975), 448--491.

\refn{[BS]}
Berele, A., and Saltman, D.J., 
{\it The centers of generic division algebras with involution}, 
Israel J. Math. (1) {\bf 63} (1988). 

\refn{[D]}
DeMeyer, F.R., 
{\it Projective modules over central separable algebras}, 
Canad. J. Math. {\bf 21} (1969), 39--43.

\refn{[K-T]}
Knus, M.A., Merkurjev, A., Rost, M., and Tignol, J.-P., 
{\it The Book of Involutions}, 
Amer. Math. Soc., Providence, RI, 1998.  

\refn{[LLT]}
Lam, T,Y., Leep, D.B., and Tignol, J.-P., 
{\it Biquaternion algebras and quartic extensions}, 
Pub. Math. IHES {\bf 77} (1993), 63--102.

\refn{[LN]}
Saltman, D.J., 
{\it Lectures on Division Algebras}, 
Amer. Math. Soc., Providence, RI, 1999.

\refn{[P]}
Procesi, C., 
{\it The invariant theory of $n \times n$ matrices}, 
Advances in Math. {\bf 19} (1976), no. 3, 306--381. 

\refn{[OS]}
Orzech, S., and Small, C., 
{\it The Brauer group of commutative rings}, 
Marcel Dekker, New York, 1975.

\refn{[R]}
Rowen, L.H., 
{\it Central simple algebras}, 
Israel J. Math. {\bf 29} (1978), 285--301. 

\refn{[R1]}
Rowen, L.H., 
{\it Polynomial identities in ring theory}, 
Academic Press, New York, London, 1980.  

\refn{[S]}
Saltman, D.J., 
{\it Retract rational fields and cyclic Galois extensions},
Israel J. Math. (2-3) {\bf 47} (1984), 165--215.

\refn{[S1]}
Saltman, D.J.,
{\it Indecomposable division algebras},
Comm. Algebra, (8) {\bf 7} (1979), 791--817.

\refn{[S2]}
Saltman, D.J., 
{\it Norm polynomials and algebras}, 
J. of Algebra, {\bf 62} (1978), 333--345.

\refn{[S3]}
Saltman, D.J., 
{\it Invariant fields of symplectic and orthogonal groups}, 
preprint

\refn{[T]}
Tate, J., 
{\it Relations between $K_2$ and Galois Cohomology}, 
Invent. Math. {\bf 36} (1976) 257--274.

\refn{[V]}
Van den Bergh, M., 
{\it The Brauer-Severi scheme of the trace ring of generic matrices}, 
Perspectives in ring theory (Antwerp 1987), 333--338;
Nato Adv. Sci. Inst. Ser. C: Math. Phys.Sci., 233, 
Kluwer Acad. Publ., Dordrecht, 1988.

\end